\documentclass{amsart}
\usepackage{amsmath} 
\usepackage{amssymb}

\newcommand{\Ga}{{\mathfrak a}}
\newcommand{\Gb}{{\mathfrak b}}
\newcommand{\Gc}{{\mathfrak c}}
\newcommand{\Gd}{{\mathfrak d}}
\newcommand{\pcf}{{\rm pcf}}
\newcommand{\uhr}{\restriction}
\newcommand{\tcf}{{\rm tcf}}

\newtheorem{theorem}{Theorem}
\newtheorem{lemma}[theorem]{Lemma} 
\newtheorem{claim}[theorem]{Proposition}

\newtheorem{subfact}{Subfact}[theorem] 

\theoremstyle{definition}
\newtheorem{definition}[theorem]{Definition}

\theoremstyle{remark}

\title{Also quite large $\mathfrak{b}\subseteq \pcf(\mathfrak{a})$ behave
nicely}
\author{Saharon Shelah}
\address{Institute of Mathematics\\
 The Hebrew University of Jerusalem\\
 91904 Jerusalem, Israel\\
 and  Department of Mathematics\\
 Rutgers University\\
 New Brunswick, NJ 08854, USA}
\email{shelah@math.huji.ac.il}
\urladdr{http://www.math.rutgers.edu/$\sim$shelah}

\begin{document}
\begin{abstract}
The present note is an answer to complains of E.Weitz on \cite{Sh:371}. We
present a corrected version of a part of chapter VIII of {\em Cardinal
Arithmetic}. 
\end{abstract}

\maketitle

\begin{definition}[{\cite[VIII 3.1]{Sh:g}}]
\label{3.1}
\begin{enumerate}
\item $J_{\ast}[\Ga]=\bigl\{\Gb:\Gb\subseteq\Gb$ and for every
inaccessible $\mu$, we have $\mu>\sup(\Gb\cap\mu)\bigr\}$.
\item $\pcf_{\ast}(\Ga)=\bigl\{\tcf(\prod\Ga/D):D$ is an ultrafilter
on $\Ga$, $D\cap J_{\ast}[\Ga]\neq\emptyset\bigr\}$. 
\item If $|\Ga|<\min(\Ga)$, for $\mu\in\pcf(\Ga)$ let $\Gb^\Ga_{\mu}=
\Gb_{\mu} [\Ga]$ be a subset of $\Ga$ such that $J_{\leq\mu}[\Ga]=J_{<\mu
}[\Ga]+\Gb_{\mu}[\Ga]$.\\
(Note that $\Gb^\Ga_{\mu}$ exists by \cite[VIII 2.6]{Sh:g}, also $\Ga$
is a finite union of $\Gb_{\mu}[\Ga]$'s).
\item If $|\Ga|<\min(\Ga)$ let $J^{\pcf}_{<\lambda}[\Ga]$ be the ideal
of subsets of $\pcf(\Ga)$ generated by $\{\pcf(\Gb_{\mu}[\Ga]):\mu\in
\lambda \cap\pcf(\Ga)\}$. 

Let $J^{\pcf}_{\leq\lambda}[\Ga]=J^{\pcf}_{<\lambda^+}[\Ga]$.
\end{enumerate}
\end{definition}

\begin{claim}[{\cite[VIII 3.1A]{Sh:g}}]
\label{3.1A}
\begin{enumerate}
\item The ideal $J^{\pcf}_{<\lambda}[\Ga]$ depends on $\Ga$ and
$\lambda$  only (and not on the choice of the $\Gb_{\mu}[\Ga]$'s). 
\item If $\Gb\subseteq\Ga$ then $J^{\pcf}_{<\lambda}[\Gb]={\mathcal
P}(\Gb)\cap J^{\pcf}_{<\lambda}[\Ga]$ and $J_{\ast}[\Gb]={\mathcal
P}(\Gb)\cap  J_{\ast}[\Ga]$. 
\end{enumerate}
\end{claim}

\begin{proof}  
(1)\quad  Let $\langle\Gb'_{\mu}[\Ga]:\mu\in\pcf(\Ga)\rangle$,
$\langle \Gb''_\mu[\Ga]:\mu\in\pcf(\Ga)\rangle$ both be as in
\ref{3.1}(3). So for each $\theta $, $\Gb'_{\theta}[\Ga]\subseteq
\Gb''_{\theta }[\Ga]\cup\bigcup_{\ell<n}\Gb''_{\theta_\ell}[\Ga]$ for
some $n<\omega$, $\theta_0,\ldots,\theta_n-1<\theta$. Hence, if
$\theta<\lambda$, 
$$
\pcf(\Gb'_{\theta}[\Ga])\subseteq\pcf(\Gb''_{\theta}[\Ga])\cup
\bigcup_{\ell< n}\pcf(\Gb''_{\theta_\ell}[\Ga]), 
$$
and each is in $J^{\pcf}_{<\lambda}[\Ga]$ as defined by
$\langle\Gb''_{\sigma}[\Ga]:\sigma\in\pcf(\Ga)\rangle$ (as
$\theta_{\ell} <\theta<\lambda$). As this\ holds for every $\theta <
\lambda$, all generators of $J^{\pcf}_{<\lambda}[\Ga]$ as defined by
$\langle \Gb'_{\sigma}[\Ga]:\sigma\in\pcf(\Ga)\rangle $ are in
$J^{\pcf}_{<\lambda}[\Ga]$ as defined by $\langle\Gb''_{\sigma}[\Ga]:
\sigma\in \pcf(\Ga)\rangle$. As the situation is symmetric we finish. 

\noindent (2) Similar proof. The first phrase follows from part (1),
and check the second. 
\end{proof}

\begin{lemma}[{\cite[VIII 3.2]{Sh:g}}]
\label{3.2} 
Suppose $|\Ga|^+<\min(\Ga)$, $\Ga\subseteq\Gb\in J_{\ast}[\pcf(\Ga)]$,
$\Gb\notin J=:J^{\pcf}_{<\lambda}[\Ga]$ and $\lambda=\max\pcf(\Ga)$. 

\noindent{\em Then\/} $\tcf(\prod\Gb/J)$ is $\lambda$.
\end{lemma}

\begin{proof}
Remember that (by \cite[VIII 2.6]{Sh:g}) there is $\langle\Gb_\theta[\Ga]:
\theta\in\pcf(\Ga)\rangle$, a generating sequence for $\Ga$. Let for
$\mu\in\pcf(\Ga)$, $\langle f^\mu_{\alpha}:\alpha<\mu\rangle$ exemplify $\mu
=\tcf(\prod\Gb_\mu [\Ga],J_{<\mu}[\Ga])$, $f^\mu _{\alpha}\in\prod\Ga$; by
\cite[3.1]{Sh:355}, without loss of generality 
\begin{enumerate}
\item[$(\ast)_0$] $\forall f\in\prod\Ga[\bigvee_{\alpha} f\uhr\Gb_\mu[\Ga]
\leq f^\mu_\alpha]$. 
\end{enumerate}
Without loss of generality for $\theta\in\Ga: f^\theta_{\alpha}(\theta)=
\alpha$ if $\alpha<\theta$, $f^\theta_{\alpha}(\theta^1)=0$ if $\alpha<
\theta<\theta^1\in\Ga$. We define $f^{\lambda,\Gb}_{\alpha}\in\prod\Gb$ by:
$$
f^{\lambda,\Gb}_{\alpha}\uhr\Ga=f^\lambda _{\alpha},
$$
and for $\theta\in\Gb\backslash\Ga$:
$$
f^{\lambda,\Gb}_{\alpha}(\theta)=\min\left\lbrace\beta: f^\lambda_\alpha
\uhr \Gb_{\theta}[\Ga]\leq f^\theta_{\beta}\mod J_{<\theta}[\Ga]
\right\rbrace. 
$$ 
Clearly 
\begin{enumerate}
\item[$(\ast)_1$] $f^\lambda_{\alpha}\leq f^\lambda_{\beta}\quad\Rightarrow
\quad f^{\lambda,\Gb}_{\alpha}\leq f^{\lambda,\Gb}_{\beta}$.
\end{enumerate}

\begin{subfact}[{\cite[VIII 3.2A]{Sh:g}}]
\label{3.2A}
$$
\alpha<\beta<\lambda\quad\Rightarrow\quad f^{\lambda,\Gb}_{\alpha}\leq
f^{\lambda,\Gb}_{\beta}\ \mod\ J.
$$
\end{subfact}

\begin{proof}[Proof of the subfact]
Let $\Gc=\{\theta\in\Ga:f^\lambda_\alpha(\theta) f^\lambda_\beta(\theta)\}$,
so $\Gc\in J_{<\lambda}[\Ga]$ and hence for some $n<\omega$ and $\sigma_1<
\ldots<\sigma_n$ from $\lambda\cap\pcf(\Ga)$ (hence $<\lambda$), we have
$\Gc \subseteq\bigcup\limits_{\ell=1}^n\Gb_{\sigma_{\ell}}[\Ga]$. So by the
definition of the $f^{\lambda,\Gb}_\alpha$'s we have: 
\begin{enumerate}
\item[$(\ast)_2$] if $\mu\in\Gb$, and $\Gb_{\mu}[\Ga]\cap\bigcap\limits_{
\ell=1}^n\Gb_{\sigma_\ell}[\Ga]\in J_{<\mu}[\Ga]$, then $f^{\lambda,
\Gb}_{\alpha}(\mu)\leq f^{\lambda,\Gb}_{\beta}(\mu)$.
\end{enumerate}
However, 
\begin{enumerate}
\item[$(\ast)_3$] $\Gd=:\left\lbrace\mu\in\pcf(\Ga):\Gb_{\mu}[\Ga]\cap
\bigcup\limits_{\ell=1}^n \Gb_{\sigma_{\ell}}[\Ga]\neq\emptyset\right.\
\mod\ J_{<\mu}[\Ga]\big\rbrace$\\
(for our fixed $\sigma_{1},\ldots,\sigma_{n}\in\lambda\cap\pcf(\Ga)$)
belongs to $J$ 
\end{enumerate}
[as $\mu\in\Gd$ implies $\mu\in\bigcup\limits_{\ell=1}^n\pcf(\Gb_{
\sigma_{\ell}}[\Ga])$ which is in $J$]. 

Together we get the subfact \ref{3.2A}.
\end{proof}

\begin{subfact}[{\cite[VIII 3.2B]{Sh:g}}]
\label{3.2B}
For any $f\in\prod\Gb$ for some $\alpha$, $f\leq f^{\lambda,\Gb}_{\alpha}$.
\end{subfact}

\begin{proof}[Proof of the subfact]
The family $J_{1}$ of sets $\Gc\subseteq\Gb$ for which this holds (i.e., for
each $f\in\prod\Gc$ there is $\alpha<\lambda$ such that $f\leq f^{\lambda,
\Gb}_{\alpha})$ satisfies: 
\begin{enumerate}
\item $\{\theta\}\in J_{1}$ for $\theta\in\lambda\cap\pcf(\Ga)$,
\item $J_{1}$ is an ideal of subsets of $\pcf(\Ga)$,
\item if $\Gc_{i}$ (for $i<\kappa$) are in $J_{1}$,  $\min(\Gc_{i})>
\kappa^+$ for each $i$ then $\bigcup\limits_{i<\kappa}\Gc_i$ is in $J_1$. 
\end{enumerate}

We shall show their satisfaction below.

This suffices for \ref{3.2B} [as $\Gb\in J_{\ast}[\pcf(\Ga)]$; why? just
prove that 
$$
\Gc\subseteq \Gb\ \&\ \Gc\in J_{\ast}[\pcf(\Ga)]\quad \Rightarrow\quad
\Gc\in J_{1} 
$$
by induction on $\sup\{\mu^+:\mu\in\Gc\}$. For successor use $(1)+(2)$. For
singular, let $\langle\mu_i:i<\kappa\rangle$ be such that $\mu_i$ is
strictly increasing continuous with limit $\sup\Gc=\sup\{\mu^+:\mu\in\Gc\}$,
and $\kappa^+<\mu_0$; by the induction hypothesis $\Gc\cap\mu_0$, $\Gc\cap
[\mu_i,\mu_{i+1}]$ are in the ideal, by (3) we know that 
$$
\bigcup_{{i}<{\kappa}}\left(\Gc\cap [\mu_{i},\mu_{i+1})\right)=\Gc\cap [\mu
_{0},\sup\Gc) 
$$
is in the ideal and by the induction hypothesis $\Gc\cap\mu_0\in J_1$ so by
(2) 
$$
\Gc=(\Gc\cap \mu _{0})\cup\left(\Gc\cap [\mu_{0},\sup\Gc)\right) 
$$
is $J_{1}$; note $\sup\Gc\notin\Gc$ as $\sup\Gc$ is singular. As $\Gb\in
 J_{\ast}[\Ga]$, we have covered all cases]. 

Now why 1), 2), 3) holds? We shall use $(\ast)_1$ from above freely.
\medskip

For (1): if $\theta\in\Ga$ as $f^{\lambda,\Gb}_\alpha\uhr\Ga=
f^\lambda_\alpha$ and $(*)_0$; if $\theta\in\Gb\setminus\Ga(\subseteq
\pcf(\Ga)$), $\alpha<\theta$ then for some $\beta<\lambda$,
$f^\theta_{\alpha+1}\leq f^\lambda_{\beta}$, hence
$f^{\lambda,\Gb}_{\beta}(\theta)>\alpha$; this shows $\{\theta\}\in J_{1}$. 

For (2): (trivially $\Gc\subseteq\Gc'\in J_1\ \Rightarrow\ \Gc\in J_1$;) if
$\Gc_{1},\Gc_{2}\in J_{1}$, $\Gc=\Gc_{1}\cup\Gc_{2}$ and $f\in\prod\Gc$,
choose, for $\ell=1,2$, $\alpha_{\ell}<\lambda$ such that $f\uhr\Gc_{\ell}
\leq f^{\lambda,\Gb}_{\alpha_{\ell}}$. Now let $f'\in\prod \Ga$ be defined
by $f'(\theta)=\max\left\lbrace f^\lambda_{\alpha_{1}}(\theta),
f^\lambda_{\alpha_{2}}(\theta)\right\rbrace$, so by an assumption on
$\langle f^\lambda_{\alpha}:\alpha<\lambda\rangle$ and $(\ast)_{0}$, for
some $\alpha$, $f'\leq f^\lambda_{\alpha}$, now $f^{\lambda,\Gb}_{\alpha}$
is as required by $(\ast)_{1}$. 

For (3): let $f\in\prod\Gc$; by assumption for each $i<\kappa$ for some
$\alpha(i)<\lambda$, $f\uhr\Gc_{i}\leq f^{\lambda,\Gb}_{\alpha(i)}$. Now
$\left(\prod\Ga,<_{J_{\leq\kappa}[\Ga]}\right)$ is $\kappa^+$-directed,
hence for some $f'\in\prod\Ga$, $\bigwedge_{{i}<{\kappa}} f^\lambda_{\alpha
(i)}<_{J_{\leq\kappa}[\Ga]}f'$. By $(\ast)_{0}$ for some $\beta<\lambda$,
$f'\leq f^\lambda_{\beta}$ and $f\uhr(\Ga\cap\Gc)\le f^\lambda_\beta$
(necessarily $ \bigwedge_{{i}<{\kappa}}\alpha(i)<\beta)$. Now for each
$\theta\in \bigcup\limits_{i<\kappa}\Gc_i$; if $\theta\in\Ga$ trivially
$f(\theta)\le f^\lambda_\beta(\theta)$, so assume $\theta\notin\Ga$; now for
some $i$, $\theta\in\Gc_i$; so $\theta>\kappa$ and $f^\lambda_{\alpha(i)}
<_{J_{\leq\kappa}[\Ga]}f^\lambda_{\beta}$, hence $f^\lambda_{\alpha(i)}
<_{J_{<\theta}[\Ga]} f^\lambda_\beta$, hence by their definitions
$f^{\lambda, \Gb}_{\alpha(i)}(\theta)\leq f^{\lambda,\Gb}_{\beta}(\theta)$. 

So $\beta$ is as required, i.e. we have proved subfact \ref{3.2B}.
\end{proof}

Now \ref{3.2} follows from \ref{3.2A}, \ref{3.2B} [using \ref{3.2B} for
$f+1$ we can get there $f<f^{\lambda,\Gb}_{\alpha}$, so (by \ref{3.2A}) for
some club $C$ of $\lambda$, 
$$
\alpha<\beta\in C\quad \Rightarrow\quad f^{\lambda,\Gb}_{\alpha}<f^{\lambda,
\Gb}_{\beta}\ \mod\ J.
$$
Together $\langle f^{\lambda,\Gb}_{\alpha}:\alpha \in C\rangle$ exemplify 
$\tcf(\prod\Gb,{<_{J}})$ is $\lambda$, as required].
\end{proof}

\begin{theorem}[{\cite[VIII 3.3]{Sh:g}}]
\label{3.3}
Assume $\min(\Ga)>|\Ga|$.
\begin{enumerate}
\item For an ultrafilter $D$ on $\pcf(\Ga)$ not disjoint to
$J_{\ast}[\pcf(\Ga)]$,  
$$
\begin{array}{ll}
\tcf(\prod\pcf(\Ga)/D) &=\min\left\lbrace\lambda\in\pcf(\Ga):
\pcf(\Gb_{\lambda}[\Ga])\in D\right\rbrace\\ 
&=\min\left\lbrace\lambda\in\pcf(\Ga):D\cap J^{\pcf}_{\leq\lambda }[\Ga] 
\neq\emptyset\right\rbrace.
\end{array}
$$
\item For $\Gc\in J_{\ast}[\pcf(\Ga)]$, $\pcf(\Gc)$ is a subset of
$\pcf(\Ga)$ and has a maximal element.
\item For $\Gb\in J_{\ast}[\pcf(\Ga)]$, $\prod\Gb/J^{\pcf}_{<\lambda}[\Ga]$
is $\lambda$-directed.
\item $\pcf_{\ast}(\Ga)=\pcf(\Ga)=\pcf_{\ast}(\pcf(\Ga))$.
\item If $\Gc\in J_{\ast}[\pcf(\Ga)]$ and $\Gc\in J^{\pcf}_{\leq\lambda}
[\Ga]$ then $\prod\Gc$ has cofinality $\leq\lambda$. 
\item If $\Gc\in J_{\ast}[\pcf(\Ga)]$ and $\Gc\in J^{\pcf}_{\leq\lambda}[\Ga] 
\backslash J^{\pcf}_{<\lambda}[\Ga]$ then $\lambda=\tcf\left(\prod\Gc,
<_{J^{\pcf}_{<\lambda }}\right)$. 
\end{enumerate}
\end{theorem}

\begin{proof} 
(1)\quad Trivially the second and third terms are equal (see Definition
\ref{3.1}(4)). Let $\lambda$ be defined as in the second term, so
$\pcf(\Gb_{\lambda}[\Ga])\in D\cap J^{\pcf}_{\leq\lambda}[\Ga]$. So by
\ref{3.1A}(2) without loss of generality $\Ga=\Gb_{\lambda}[\Ga]$, so
$\lambda=\max\pcf(\Ga)$. Using \ref{3.2}'s notation, $\langle f^{\lambda,
\Gb}_{\alpha}: \alpha<\lambda\rangle$ exemplify $\lambda=\tcf(\prod\Ga/D)$.

\noindent (2)\quad By (1).

\noindent (3)\quad This follows by the proof of \ref{3.2}, but as I was
asked, we repeat the proof of \ref{3.2} with the required
changes. W.l.o.g. $\lambda\in\pcf(\Ga)$ [why? if $\lambda>\max\pcf(\Ga)$
then $J^{\pcf}_{<\lambda}[\Ga]={\mathcal P}(\pcf(\Ga))$, so the conclusion is
trivial, if not let $\lambda'=\min(\pcf(\Ga)\setminus\lambda)$, so $\lambda'
\in\pcf(\Ga)$ and $J^{\pcf}_{<\lambda}[\Ga]=J^{\pcf}_{<\lambda'}[\Ga]$]. We
let $J=:J^{\pcf}_{<\lambda}[\Ga]$. Remember that (by \cite[ VIII 2.6]{Sh:g})
there is $\langle\Gb_\theta[\Ga]:\theta\in\pcf(\Ga)\rangle $, a generating
sequence for $\Ga$. Let for $\mu\in\pcf(\Ga)$, $\langle f^\mu _{\alpha}:
\alpha<\mu\rangle$ exemplify $\mu=\tcf(\prod\Gb_\mu[\Ga],J_{<\mu}[\Ga])$,
$f^\mu_{\alpha}\in\prod\Ga$; by\cite[3.1]{Sh:355} without loss of generality
\begin{enumerate}
\item[$(\ast)_{0}$] $\forall f\in \prod \Ga[\bigvee_{\alpha} f\uhr\Gb_\mu
[\Ga] \leq f^\mu _{\alpha}]$. 
\end{enumerate}
Without loss of generality for $\theta\in\Ga: f^\theta_{\alpha}(\theta)=
\alpha$ if $\alpha<\theta$, $f^\theta_\alpha(\theta^1)=0$ if $\alpha<\theta
<\theta^1\in\Ga$. For any $f\in\prod\Ga$ we define a function $f^{\Gb}\in
\prod \Gb$ by: 
$$
f^{\Gb}\uhr\Ga=f
$$
and for $\theta\in\Gb\backslash\Ga$:
$$
f^{\Gb}(\theta)=\min\left\lbrace\beta:f\uhr \Gb_{\theta}[\Ga]\leq
f^\theta_{\beta}\; \mod\; J_{<\theta}[\Ga]\right\rbrace.
$$
Let $f$ vary on $\prod \Ga$. Clearly 
\begin{enumerate}
\item[$(\ast)_1$] $f_1\leq f_2\quad\Rightarrow\quad f^\Gb_1\leq f^\Gb_2$. 
\end{enumerate}

\begin{subfact}
\label{3.3A}
If $f_1\leq f_2\ \mod\ J_{<\lambda}[\Ga]$ (both in $\prod\Ga$ of course)
then $f^\Gb_1\leq f^\Gb_2\ \mod\ J$. 
\end{subfact}

\begin{proof}[Proof of the subfact]
Let $\Gc=\{\theta\in\Ga:f_1(\theta)\geq f_2(\theta)\}$, so $\Gc\in
J_{<\lambda}[\Ga]$, hence for some $n<\omega$ and $\sigma_1<\ldots<
\sigma_n$ from $\lambda\cap\pcf(\Ga)$, (hence $<\lambda$) we have $\Gc
\subseteq\bigcup\limits_{\ell=1}^n \Gb_{\sigma_{\ell}}[\Ga]$. So by the
definition of the $f^{\Gb}_\ell$'s we have: 
\begin{enumerate}
\item[$(\ast)_{2}$] if $\mu\in\Gb$, and $\Gb_{\mu}[\Ga]\cap
\bigcap\limits_{\ell=1}^n\Gb_{\sigma_{\ell}}[\Ga]\in J_{<\mu}[\Ga]$, then
$f^{\Gb}_{1}(\mu)\leq f^{\Gb}_{2}(\mu)$.
\end{enumerate}
However 
\begin{enumerate}
\item[$(\ast)_{3}$] $\Gd=:\left\lbrace\mu\in \pcf(\Ga):\Gb_{\mu}[\Ga]\cap
\bigcup\limits_{\ell=1}^n \Gb_{\sigma_{\ell}}[\Ga]\neq\emptyset\ \mod\
J_{<\mu }[\Ga]\right\rbrace$\\  
(for our fixed $\sigma_{1},\ldots,\sigma_{n}\in\lambda\cap\pcf(\Ga)$)
belongs to $J$
\end{enumerate}
[as $\mu\in\Gd$ implies $\mu\in\pcf(\bigcup\limits^n_{\ell=1}
\Gb_{\sigma_\ell}[\Ga])=\bigcup\limits_{\ell=1}^n\pcf(\Gb_{\sigma_{\ell}}
[\Ga])$ which is in $J$]. 

Together we get subfact \ref{3.3A}.
\end{proof}

\begin{subfact}
\label{3.3B}
For any $g\in\prod(\Gb)$ for some $f\in\prod(\Ga)$ we have $g\leq
f^{\lambda, \Gb}_{\alpha}$.
\end{subfact}

\begin{proof}[Proof of the subfact]
The family $J_{1}$ of sets $\Gc\subseteq\Gb$ for which this holds, i.e., for
each $g\in\prod\Gc$ there is $f\in\prod\Ga$ such that $g\leq f^{\Gb}$,
satisfies: 
\begin{enumerate}
\item $\{\theta\}\in J_{1}$ for $\theta\in\lambda\cap\pcf(\Ga)$,
\item $J_{1}$ is an ideal of subsets of $\pcf(\Ga)$,
\item if $\Gc_i$ (for $i<\kappa$) are in $J_1$, $\min(\Gc_i)>\kappa^+$ for
each $i$ {\em then\/} $\bigcup\limits_{i<\kappa}\Gc_i$ is in $J_1$.
\end{enumerate}

We shall show their satisfaction below.

{\em Why (1)+(2)+(3) suffice for \ref{3.3B}?}\quad As $\Gb\in
J_{\ast}[\pcf(\Ga)]$; why? just prove that
\begin{enumerate}
\item[$(*)_4$] $\Gc\subseteq\Gb\ \&\ \Gc\in J_{\ast}[\pcf(\Ga)]\quad
\Rightarrow\quad\Gc\in J_{1}$
\end{enumerate}
by induction on $\sup\{\mu^+:\mu\in\Gc\}$. For successor use $(1)+(2)$. For
singular, let $\langle\mu_i:i<\kappa\rangle$ be such that $\mu_i$ is
strictly increasing continuous with limit $\sup(\Gc)=\sup\{\mu^+:\mu\in\Gc
\}$, and $\kappa^+<\mu_0$; by the induction hypothesis $\Gc\cap\mu_0$,
$\Gc\cap [\mu_{i},\mu_{i+1}]$ are in the ideal, by (3) we know that 
$$
\bigcup_{{i}<{\kappa}}\left(\Gc\cap [\mu_{i},\mu_{i+1})\right)=\Gc\cap
[\mu_{0},\sup\Gc) 
$$
is in the ideal and, as said above, $\Gc\cap\mu_0\in J_1$ so by (2)
$$
\Gc=(\Gc\cap\mu_{0})\cup\left(\Gc\cap[\mu_{0},\sup\Gc)\right)
$$
is $J_{1}$; note $\sup(\Gc)\notin\Gc$ as $\sup(\Gc)$ is singular. As $\Gc\in
J_{\ast}[\pcf(\Ga)]$ implies $\Gc$ has no inaccessible accumulation point,
we have covered all cases in the induction, so $(*)_4$ holds. Now note that
$\Gb\in J_*[\pcf(\Ga)]$, so from $(*)_4$ we get $\Gb\in J_1$ and by the
definition of $J_1$ we are done. 
\medskip

{\em Next, why 1), 2), 3) hold?}\quad We shall use $(\ast)_1$ from above
freely.  

For (1): let $g\in\prod\Gb$; if $\theta\in\Ga$ as $f^{\Gb}\uhr \Ga=f$ and
$(*)_0$; if $\theta\in\Gb\setminus\Ga$ ($\subseteq\pcf(\Ga)$), then
$g(\theta)<\theta$ let $f=f^\theta_{g(\theta)+1}$ hence 
$$
(\forall\gamma\leq g(0))(\neg f\leq f^\theta_\gamma\ \mod\
J_{<\theta}[\Ga]). 
$$
Hence $g(\theta)<f^{\Gb}(\theta)$; this shows $\{\theta\}\in J_{1}$. 

For (2): (trivially $\Gc\subseteq \Gc'\in J_1\ \Rightarrow\ \Gc\in J_1$); if
$\Gc_1,\Gc_2\in J_{1}$, $\Gc=\Gc_{1}\cup\Gc_{2}$ and $g\in\prod\Gc$, choose
for $\ell=1,2$ $f_{\ell}\in\prod\Ga$ such that $g\uhr\Gc_{\ell}\leq
f^{\Gb}_{\ell}$. Now let $f\in\prod\Ga$ be defined by $f(\theta)=\max
\left\lbrace f_{1}(\theta), f_{2}(\theta)\right\rbrace$, so $f\in\prod\Ga$
and $g\uhr\Gc_1\leq f^\Gb_1\leq f^\Gb$ and $g\uhr\Gc_2\leq f^\Gb_2\leq
f^\Gb$ hence $g\uhr (\Gc_1\cup \Gc_2)\leq f^\Gb$. 

For (3): let $g\in\prod\Gc$; by assumption for each $i<\kappa $ for some
$f_i\in\prod\Ga$, $g\uhr\Gc_{i}\leq f^{\Gb}_{i}$. Now $\left(\prod\Ga,
<_{J_{\leq\kappa}[\Ga]}\right)$ is $\kappa^+$-directed, hence for some $f\in
\prod\Ga$, $\bigwedge_{{i}<{\kappa}} f_{i}<_{J_{\leq\kappa}[\Ga]}f$.
W.l.o.g. $g\uhr\Ga\leq f\uhr\Gc$. Now for each $\theta\in
\bigcup\limits_{i<\kappa}\Gc_i$; if $\theta\in\Ga$ trivially $g(\theta)\leq
f(\theta)$, so assume $\theta\notin\Ga$. Now, for some $i$, $\theta\in
\Gc_{i}$; so $\theta>\kappa$ and $f_i<_{J_{\leq\kappa}[\Ga]}f$, hence $f_{i}
<_{J_{<\theta}[\Ga]} f$, hence by their definitions $f^{\Gb}_{i}(\theta)
\leq f^{\Gb}(\theta)$. 
\medskip

So (1), (2), (3) hold and hence $\Gb$ is as required, i.e., we have proved
subfact \ref{3.3B}.  
\end{proof}

We finish by

\begin{subfact}
\label{3.3C}
$\prod \Gb/J$ is $\lambda$-directed.
\end{subfact}

\begin{proof}[Proof of the subfact]
Assume $g_i\in\prod\Gb$ for $i<i^*<\lambda$, by \ref{3.3B} for each $i<i^*$
for some $f_i\in\prod\Ga$ we have $g_i\leq f^{\Gb}_i$. But
$\prod\Ga/J_{<\lambda}[\Ga]$ is $\lambda$-directed hence for some $f\in
\prod\Ga$ we have  
$$
\bigwedge_{i<i^*} f_i<f\ \mod\ J_{<\lambda}[\Ga].
$$
By \ref{3.3A} we have 
$$
\bigwedge_{i<i^*} f_i^\Gb\leq f^\Gb\ \mod\ J,
$$
hence by the previous sentence $i<i^*\quad \Rightarrow\quad g_i\leq f^\Gb_i
\leq_J f^\Gb$, so $f^\Gb+1$ is a $<_J$-upper bound of $\{g_i:i<i^*\}$, as
required. 
\end{proof}
\end{proof}

%\bibliographystyle{literal-unsrt}
%\bibliography{lista,listb,listx}

\begin{thebibliography}{Sh 371}
\makeatletter \renewcommand{\@biblabel}[1]{[#1]} \makeatother

\bibitem[Sh 371]{Sh:371}Saharon Shelah.
\newblock {Advanced: cofinalities of small reduced products}.
\newblock In {\em {Cardinal Arithmetic}}, volume~29 of {\em {Oxford Logic
  Guides}}, chapter {VIII}. {Oxford University Press}, 1994.

\bibitem[Sh:g]{Sh:g}Saharon Shelah.
\newblock {\em {Cardinal Arithmetic}}, volume~29 of {\em {Oxford Logic
  Guides}}.
\newblock {Oxford University Press}, 1994.

\bibitem[Sh 355]{Sh:355}Saharon Shelah.
\newblock {$\aleph _{\omega +1}$ has a Jonsson Algebra}.
\newblock In {\em {Cardinal Arithmetic}}, volume~29 of {\em {Oxford Logic
  Guides}}, chapter~II. {Oxford University Press}, 1994.

\end{thebibliography}

\end{document}